\documentclass[preprint]{imsart}

\usepackage{amsthm,amsmath,amssymb,amsbsy}
\usepackage[authoryear]{natbib}
\startlocaldefs
	\newtheorem{thm}{Theorem}[section]
	\newtheorem{cor}[thm]{Corollary}
	\newtheorem{lem}[thm]{Lemma}
	
	\newtheorem{cond}[thm]{Condition}
	
	\theoremstyle{definition}
	
	\newtheorem{rem}[thm]{Remark}

    \newcommand{\norm}[1]{\| #1 \|}
    \newcommand{\eps}{\varepsilon}
    
    \newcommand{\RR}{\mathbb{R}}
    \newcommand{\Sb}{\mathbb{S}}
    \newcommand{\EE}{\mathbb{E}}
    \newcommand{\NN}{\mathbb{N}}
    \newcommand{\ZZ}{\mathbb{Z}}

    \newcommand{\law}{\mathcal{L}}

    \renewcommand{\a}{{\boldsymbol{a}}}
    \newcommand{\A}{\boldsymbol{A}}
    \newcommand{\B}{\boldsymbol{B}}
    \newcommand{\C}{\boldsymbol{C}}
    \newcommand{\Q}{\boldsymbol{Q}}
    \newcommand{\X}{\boldsymbol{X}}
    \newcommand{\Y}{\boldsymbol{Y}}
    \newcommand{\Z}{\boldsymbol{Z}}
    
    \newcommand{\cc}{\boldsymbol{c}}
    \newcommand{\bfe}{\boldsymbol{e}}

    \newcommand{\x}{\boldsymbol{x}}
    \newcommand{\y}{\boldsymbol{y}}

    \newcommand{\0}{\boldsymbol{0}}
    \newcommand{\1}{\boldsymbol{1}}
    \newcommand{\xib}{\boldsymbol{\xi}}

    \newcommand{\E}{\mathrm{E}}

    \newcommand{\Thetab}{\boldsymbol{\Theta}}

    \newcommand{\dto}{\rightsquigarrow}%{\stackrel{d}{\to}}
    \newcommand{\wto}{\stackrel{w}{\to}}
    \newcommand{\vto}{\stackrel{v}{\to}}
    \renewcommand{\d}{\mathrm{d}}

    \newcommand{\rmd}{\,\mathrm{d}}
    
    \newcommand{\dy}{\rmd y}

    \newcommand{\be}{\begin{equation}}
    \newcommand{\ee}{\end{equation}}
	
	\numberwithin{equation}{section}
\endlocaldefs

\begin{document}

\begin{frontmatter}

% "Title of the paper"
%\title{Regularly varying multivariate time series}
\title{Regularly varying \\ multivariate time series}
\runtitle{Multivariate time series}

\author{\fnms{Bojan} \snm{Basrak}\corref{}\ead[label=bojan]{bbasrak@math.hr}\thanksref{MZOS} and 
\fnms{Johan} \snm{Segers}\ead[label=johan]{Johan.Segers@uclouvain.be}\thanksref{IAP}}
\affiliation{University of Zagreb and Universit\'e catholique de Louvain}
\thankstext{MZOS}{Supported by research grant MZO\v S project nr.~037-0372790-2800 of the Croatian government.}
\thankstext{IAP}{Supported by a post-doctoral grant of the Fund for Scientific Research, Flanders, a VENI grant of the Netherlands Organization for Scientific Research (NWO) and by the IAP research network grant nr.\ P6/03 of the Belgian government (Belgian Science Policy).}
\address{University of Zagreb \\ Department of Mathematics \\ Bijeni\v cka 30 \\ 10000, Zagreb (Croatia) \\ \printead{bojan}}
\address{Universit\'e catholique de Louvain \\ Institut de statistique \\ Voie du Roman Pays, 20 \\ B-1348 Louvain-la-Neuve (Belgium) \\ \printead{johan}}

\runauthor{B. Basrak and J. Segers}

\begin{abstract}
A multivariate, stationary time series is said to be jointly regularly varying if all its finite-dimensional distributions are multivariate regularly varying. This property is shown to be equivalent to weak convergence of the conditional distribution of the rescaled series given that, at a fixed time instant, its distance to the origin exceeds a threshold tending to infinity. The limit object, called the tail process, admits a decomposition in independent radial and angular components. Under an appropriate mixing condition, this tail process allows for a concise and explicit description of the limit of a sequence of point processes recording both the times and the positions of the time series when it is far away from the origin. The theory is applied to multivariate moving averages of finite order with random coefficient matrices.
\end{abstract}



\end{frontmatter}

\section{Introduction}
\label{S:intro}

Extreme values of a stationary, multivariate time series may exhibit dependence across coordinates and over time. Multivariate regular variation, on the one hand, is the central concept in the theory of extremes of random vectors \citep{Resnick87, Resnick06}. Point processes, on the other hand, provide a convenient language for describing temporal dependence of time series extremes \citep{LeRo88, HsHuLe89}. This explains why in Davis and Mikosch (1998), building upon \citet{DaHs95}, regular variation and point process techniques are the two main ingredients of a theory of extremes of multivariate time series. The dynamics of the series in regions far away from the origin are reflected in the typical property that the atoms of the limiting point processes can be grouped into independent clusters. The only handle on the distribution of these clusters, however, is in the form of a number of asymptotic relations that are difficult to evaluate in general.
%\citep[Theorems~2.6 and 2.8]{DaMi98}

One of our aims is then to reconstruct and extend the theory in Davis and Hsing (1995) and Davis and Mikosch (1998) in such a way that the relation between the original time series and the limiting point process of extremes becomes explicit. The argument rests on a novel characterization of the property that the finite-dimensional distributions of a discrete-time, stationary time series $(\X_t)_{t \in \ZZ}$ in $\RR^d$ are multivariate regularly varying. The necessary and sufficient condition is that, for a given norm $\norm{\,\cdot\,}$ on $\RR^d$, the conditional distribution of the rescaled series $(x^{-1} \X_t)_{t \in \ZZ}$ given $\norm{\X_0} > x$ converges to a non-degenerate limit as $x \to \infty$. The limit process $(\Y_t)_{t \in \ZZ}$, called the tail process of $(\X_t)$, admits a decomposition in two independent components: on the one hand, a real-valued radial component $\norm{\Y_0}$ whose distribution only depends on the tail index $\alpha \in (0, \infty)$ of $\norm{\X_0}$; on the other hand, a sequence-valued angular component $(\Thetab_t)_{t \in \ZZ}$ defined by $\Thetab_t = \Y_t / \norm{\Y_0}$ for $t \in \ZZ$, called the spectral process. Moreover, the distributions of the restrictions of the spectral process $(\Thetab_t)$ to positive or negative times are connected by a certain adjoint relation parametrized by $\alpha$.

Many questions concerning excursions of $(\X_t)$ far away from the origin can be answered via a study of the sequence of time-space point processes
\begin{equation}
\label{E:Nn}
    N_n = \sum_{i=1}^n \delta_{(i/n, \X_i/a_n)},
\end{equation}
with $(a_n)$ a positive sequence such that $n \Pr(\norm{\X_0} > a_n) \to 1$ as $n \to \infty$. Under a tailor-made form of mixing condition, the atoms of the limiting point process $N$ can be partitioned into independent clusters. The distribution of these clusters can be explicitly and concisely described in terms of the tail process $(\Y_t)$ of $(\X_t)$. Among other things, these results lead to convenient formulas of the extremal index and cluster size probabilities of various univariate series derived from $(\X_t)$.

The results are applied to a multivariate moving average process of finite order with random coefficient matrices, defined for $t \in \ZZ$ by
\begin{equation}
\label{E:MMA}
    \X_t = \sum_{i=0}^m \C_i(t) \xib_{t-i} ;
\end{equation}
here $(\xib_t)_{t \in \ZZ}$ is a sequence of independent and indentically distributed random column vectors in $\RR^q$ and $\{ \C_i(t) : i = 0, \ldots, m ; t \in \ZZ \}$ is an array of random $d \times q$ matrices, independent of $(\xib_t)$ and stationary as a $d \times q \times (m+1)$ dimensional process indexed over $t \in \ZZ$. The tail behaviour of the stationary distribution of infinite-order versions of \eqref{E:MMA} has been studied in \citet{ReWi91} and \citet{HuSa07}; in contrast, our focus is on the process as a whole but under the simplifying assumption that the order is finite and $\{ \C_i(t) \}$ is independent of $(\xib_t)$.

Main results are gathered in section~\ref{S:main} and are grouped into the three themes already identified in the previous paragraphs: regular variation and the tail process (subsection~\ref{SS:main:tailprocess}), point processes of extremes (subsection~\ref{SS:main:pointprocess}), and moving averages with random coefficients (subsection~\ref{SS:main:MMA}). Proofs of the main theorems as well as statements and proofs of additional results are given in sections \ref{S:tailprocess} to \ref{S:mma} for the same three themes respectively. Additional results we wish to mention here concern Laplace functionals of point processes of clusters of extremes in Theorem~\ref{T:Laplace} and extremal indices of linear projections of multivariate time series in Remark~\ref{rem:linear}.

The central notion in the paper is that of regular variation. The law of a $d$-dimensional random vector $\X$ is called regularly varying of index $\alpha \in (0, \infty)$ if for some norm $\norm{\,\cdot\,}$ on $\RR^d$ there exists a random vector $\Thetab$ on the unit sphere $\Sb^{d-1} = \{ \x \mid \norm{\x} = 1 \}$ such that for every $u \in (0, \infty)$ and as $x \to \infty$,
\begin{equation}
\label{E:RV}
    \frac{1}{\Pr(\norm{\X} > x)} \Pr( \norm{\X} > ux, \, \X / \norm{\X} \in \cdot \,)
    \wto u^{-\alpha} \Pr(\Thetab \in \cdot \,),
\end{equation}
where $\wto$ denotes weak convergence of finite measures. The law of $\Thetab$ is called the spectral measure of $\X$. The definition of regular variation does not depend on the particular norm chosen in the sense that \eqref{E:RV} holds for some norm if and only if it holds for every norm, the spectral measure of course depending on the norm. Finally, a stationary $d$-dimensional time series $(\X_t)_{t \in \ZZ}$ is said to be jointly regularly varying of index $\alpha \in (0, \infty)$ if for every positive integer $k$ the $kd$-dimensional random vector $(\X_1, \ldots, \X_k)$ is regularly varying of index $\alpha$. Some other references on multivariate regular variation, apart from the ones already mentioned, are \citet{BDM02a}, \citet{HuLi06}, and \citet{MS01}.

Besides the notation already appearing in this introduction, the following symbols are used throughout the paper: $\NN = \ZZ \cap [0, \infty)$, $\EE = [-\infty, \infty]^d \setminus \{ \0 \}$ and $\EE_u = \{ \x \in \EE \mid \norm{\x} > u \}$; the law of the random vector $\X$ is $\law(\X)$; the indicator variable of an event $E$ is $\1(E)$; convergence of probability distributions and vague convergence of Radon measures (see section~\ref{S:tailprocess}) is indicated by $\dto$ and $\vto$, respectively; for a topological space $\mathbb{T}$ the space of continuous functions $f: \mathbb{T} \to \RR$ is denoted by $C(\mathbb{T})$, decorations with the subscript $K$ and the superscript + indicating the subclasses of those $f \in C(\mathbb{T})$ that have compact support or take values in $[0, \infty)$, respectively.

% ===============================================================================
\section{Main results}
\label{S:main}

In subsection~\ref{SS:main:tailprocess}, joint regular variation of a stationary time series $(\X_t)$ is identified with a certain asymptotic property of the conditional distribution of the series given that it is far away from the origin at a fixed time instant. This characterization is exploited in subsection~\ref{SS:main:pointprocess} to describe limits of certain point processes of extremes. An application to multivariate moving averages with random coefficient matrices is given in subsection~\ref{SS:main:MMA}. The proofs of the results in this section are to be found further on in the paper.

%%%%%%%%%%%
\subsection{Tail process}
\label{SS:main:tailprocess}

The most important object in this paper is introduced in our first theorem.

\begin{thm}
\label{T:JRV}
Let $(\X_t)_{t \in \ZZ}$ be a stationary process in $\RR^d$ and let $\alpha \in (0, \infty)$. The following three statements are equivalent:
\begin{itemize}
\item[(i)]
$(\X_t)$ is jointly regularly varying of index $\alpha$.
\item[(ii)]
There exists a process $(\Y_t)_{t \in \NN}$ in $\RR^d$ with $\Pr(\norm{\Y_0} > y) = y^{-\alpha}$ for $y \geq 1$ such that for every $t \in \NN$ and as $x \to \infty$,
\[
    \law( x^{-1} \X_0, \ldots, x^{-1} \X_t \mid \norm{\X_0} > x )
    \dto \law( \Y_0, \ldots, \Y_t ).
\]
\item[(iii)]
There exists a process $(\Y_t)_{t \in \ZZ}$ in $\RR^d$ with $\Pr(\norm{\Y_0} > y) = y^{-\alpha}$ for $y \geq 1$ such that for all $s, t \in \ZZ$ with $s \leq t$ and as $x \to \infty$,
\[
    \law( x^{-1} \X_s, \ldots, x^{-1} \X_t \mid \norm{\X_0} > x )
    \dto \law( \Y_s, \ldots, \Y_t ).
\]
\end{itemize}
\end{thm}

The process $(\Y_t)_{t \in \ZZ}$ in Theorem~\ref{T:JRV}(iii) is called the tail process of $(\X_t)$. In general, the tail process is itself not stationary. It has a number of remarkable properties, the two most important ones being described next.

\begin{thm}
\label{T:tailprocess}
Let $(\Y_t)_{t \in \ZZ}$ be the tail process in Theorem~\ref{T:JRV}(iii) and define $\Thetab_t = \Y_t / \norm{\Y_0}$ for $t \in \ZZ$.

(i) $\norm{\Y_0}$ is independent of $(\Thetab_t)_{t \in \ZZ}$.

(ii) For all $i,s,t \in \ZZ$ with $s \leq 0 \leq t$ and for all bounded and continuous  $f : (\RR^d)^{t-s+1} \to \RR$ satisfying $f(\y_s, \ldots, \y_t) = 0$ whenever $\y_0 = 0$,
\begin{equation}
\label{E:shift}
    \E [ f(\Thetab_{s-i}, \ldots, \Thetab_{t-i}) ] =
    \E  \biggl[
        f   \biggl(
                \frac{\Thetab_{s}}{\norm{\Thetab_{i}}}, \ldots,
                \frac{\Thetab_t}{\norm{\Thetab_{i}}}
            \biggr)
        \norm{\Thetab_{i}}^\alpha
        \biggr].
\end{equation}
\end{thm}

Theorem~\ref{T:tailprocess}(i) shows that the distribution of $(\Y_t)$
can be decomposed into a real-valued radial component, $\norm{\Y_0}$, and a sequence-valued angular component, $(\Thetab_t)_{t \in \ZZ}$, the two components being independent. As the law of $\Thetab_0$ is the spectral measure of the one of $\X_0$, we coin the process $(\Thetab_t)$ the spectral process of $(\X_t)$. The characterizations of joint regular variation in Theorem~\ref{T:JRV} can be rephrased in terms of this spectral process.

\begin{cor}
\label{cor:spectralprocess}
Let $(\X_t)_{t \in \ZZ}$ be a stationary process in $\RR^d$. Assume that the function $x \mapsto \Pr(\norm{\X_0} > x)$ is regularly varying of index $-\alpha$ for some $\alpha \in (0, \infty)$. The following three statements are equivalent:
\begin{itemize}
\item[(i)]
$(\X_t)$ is jointly regularly varying of index $\alpha$.
\item[(ii)]
There exists a process $(\Thetab_t)_{t \in \NN}$ in $\RR^d$ such that for every $t \in \NN$ and as $x \to \infty$,
\[
    \law \biggl( \frac{\X_0}{\norm{\X_0}}, \ldots, \frac{\X_t}{\norm{\X_0}}
        \, \bigg| \, \norm{\X_0} > x \biggr)
    \dto \law( \Thetab_0, \ldots, \Thetab_t ).
\]
\item[(iii)]
There exists a process $(\Thetab_t)_{t \in \ZZ}$ in $\RR^d$ such that for every $s, t \in \ZZ$ with $s \leq t$ and as $x \to \infty$,
\[
    \law \biggl( \frac{\X_s}{\norm{\X_0}}, \ldots, \frac{\X_t}{\norm{\X_0}}
        \, \bigg| \, \norm{\X_0} > x \biggr)
    \dto \law( \Thetab_s, \ldots, \Thetab_t ).
\]
\end{itemize}
In this case, the tail process $(\Y_t)$ of $(\X_t)$ is given by $\Y_t = Y \Thetab_t$ for $t \in \ZZ$, the random variable $Y$ being independent of $(\Thetab_t)$ and having survival function $\Pr(Y > y) = y^{-\alpha}$ for $y \in [1, \infty)$.
\end{cor}

Further, for $j \in \NN$, Theorem~\ref{T:tailprocess}(ii) with on the one hand $s = 0$ and $t = i = j$ and on the other hand $s = i = -j$ and $t = 0$ stipulates that the distributions of $(\Thetab_{-j}, \ldots, \Thetab_0)$ and $(\Thetab_0, \ldots, \Thetab_j)$ are in some sense adjoint to each other. For univariate Markov chains, this adjoint relation was already described in \citet{Segers07}. 

The proofs of Theorems~\ref{T:JRV} and \ref{T:tailprocess} and Corollary~\ref{cor:spectralprocess} are given in section~\ref{S:tailprocess}, together with some further properties of the spectral process.

%%%%%%%%%%%
\subsection{Point processes}
\label{SS:main:pointprocess}

Throughout this subsection, let $(\X_t)$ be a stationary time series in $\RR^d$, regularly varying of index $\alpha \in (0, \infty)$ and with tail process $(\Y_t)$ and spectral process $(\Thetab_t)$. Also write $M_r = \max_{i = 1, \ldots, r} \norm{\X_i}$ for $r \in \NN$. Let $(a_n)$ be a positive sequence such that $n \Pr(\norm{\X_0} > a_n) \to 1$ as $n \to \infty$. Of interest is the weak limit of the time-space point processes
\[
    N_n = \sum_{i=1}^n \delta_{(i/n, \X_i/a_n)}.
\]
In order to control the total mass on the time axis, the state space needs to be restricted to $[0,1]\times \EE_u$ for some $u \in (0, \infty)$. 

An important role will be played by the quantity $\theta$ defined by
\begin{eqnarray}
\label{E:candidatetheta}
	\lefteqn{
	\lim_{r \to \infty} \lim_{x \to \infty}
	\Pr( M_r \leq x \mid \norm{\X_0} > x )
	} \\
	&=& \lim_{r \to \infty} \Pr \biggl( \max_{i = 1, \ldots, r} \norm{\Y_i} \leq 1 \biggr) \nonumber
	= \Pr \biggl( \sup_{i \geq 1} \norm{\Y_i} \leq 1 \biggr)
	=: \theta. 
\end{eqnarray}
In view of what is to follow, $\theta$ is coined the candidate extremal index of the univariate series $(\norm{\X_t})$. By Theorem~\ref{T:tailprocess}(i),
\begin{eqnarray}
\label{E:theta:spectral}
    \theta
    &=& \int_1^\infty 
    \Pr \biggl( \sup_{i \geq 1} \norm{\Thetab_i}^\alpha \leq y^{-\alpha} \biggr)
    \rmd(-y^{-\alpha}) \\
    &=& \E \max \biggl( 1 - \sup_{i \geq 1} \norm{\Thetab_i}^\alpha, 0 \biggr)
    = \E \biggl[ \sup_{i \geq 0} \norm{\Thetab_i}^\alpha
    - \sup_{i \geq 1} \norm{\Thetab_i}^\alpha \biggr]. \nonumber
\end{eqnarray}
Further, let $E$ be the event that the supremum of $(\norm{\Y_t})$ is attained for the first time at $t = 0$, that is,
\begin{equation}
\label{E:E}
    E = \biggl\{ \sup_{t < 0} \norm{\Y_t} < \norm{\Y_0} = \sup_{t \in \ZZ} \norm{\Y_t} \biggr\}.
\end{equation}

The following condition prohibits clusters of extremes to linger on for too long \citep{DaHs95}.

\begin{cond}
\label{C:FiCl}
There exists a positive integer sequence $(r_n)$ such that $r_n \to \infty$ and $r_n / n \to 0$ as $n \to \infty$ and such that for every $u \in (0, \infty)$,
\[
    \lim_{m\to\infty} \limsup_{n\to\infty}
    \Pr \biggl( \max_{m \leq |t| \leq r_n} \norm{\X_t} > a_n u \, \bigg| \, \norm{\X_0} > a_n u \biggr) = 0.
\]
\end{cond}

The following theorem then describes the limit distribution of point processes of clusters of extremes.

\begin{thm}
\label{T:pointprocess:cluster}
If Condition~\ref{C:FiCl} holds, then, with $\theta$ and $E$ as in \eqref{E:candidatetheta} and \eqref{E:E} respectively, $\Pr( \lim_{|t| \to \infty} \norm{\Y_t} = 0 ) = 1$ and $\theta = \Pr(E) > 0$. Moreover, for every $u \in (0, \infty)$ and as $n \to \infty$,
\begin{equation}
\label{E:runsblocks}
    \Pr(M_{r_n} \leq a_n u \mid \norm{\X_0} > a_n u)
    = \frac{\Pr(M_{r_n} > a_n u)}{r_n \Pr(\norm{\X_0} > a_n u)} + o(1) 
	\to \theta
\end{equation}
and in the state space $\EE$,
\begin{equation}
\label{E:clusterprocess}
    \law \Biggl( \sum_{i=1}^{r_n} \delta_{\X_i / (a_n u)} \, \Bigg| \, M_{r_n} > a_n u \Biggr)
    \dto \law \Biggl( \sum_{j=1}^\infty \delta_{\Z_j} \Biggr),
\end{equation}
the law of $\sum_{j = 1}^\infty \delta_{\Z_j}$ being equal to the one of $\sum_{t \in \ZZ} \delta_{\Y_t}$ conditionally on $E$.
\end{thm}

\begin{rem}[mean cluster size]
\label{rem:meanclustersize}
A combination of \eqref{E:runsblocks} and \eqref{E:clusterprocess} shows that if Condition~\ref{C:FiCl} holds, then for every $u \in (0, \infty)$ and as $n \to \infty$,
\[
	\E \Biggl[ \sum_{i=1}^{r_n} \1(\norm{\X_i} > a_n u) \, \Bigg| \, M_{r_n} > a_n u \Biggr] \to \E \Biggl[ \sum_{t \in \ZZ} \1(\norm{\Y_t} > 1) \, \Bigg| \, E \Biggr] 
	= \frac{1}{\theta} < \infty.
\]
\end{rem}

In order to describe the limit of $N_n$, the following extension of condition
$\mathcal{A}(a_n)$ of \citet{DaHs95} is needed.
We note that both of them are implied
by the strong mixing property.

\begin{cond}[$\mathcal{A}'(a_n)$]
\label{C:Aprimean}
There exists a positive integer sequence $(r_n)$ such that $r_n \to \infty$ and $r_n / n \to 0$ as $n \to \infty$ and such that for every $f \in C_K^+([0,1] \times \EE)$, denoting $k_n = \lfloor n/r_n \rfloor$,
\[
    \E \exp \Biggl\{ - \sum_{i=1}^n f \biggl(\frac{i}{n},\frac{\X_i}{a_n}\biggr) \Biggr\}
    - \prod_{k=1}^{k_n} \E \exp \Biggl\{ - \sum_{i=1}^{r_n} f \biggl(\frac{k r_n}{n}, \frac{\X_i}{a_n}\biggr) \Biggr\}
    \; \,\to\ \, 0\,.
\]
\end{cond}

\begin{thm}
\label{T:pointprocess:complete}
If Conditions~\ref{C:FiCl} and \ref{C:Aprimean} hold, then for every $u \in (0, \infty)$ and as $n \to \infty$,
\[
    N_n \dto N^{(u)} 
    = \sum_i \sum_j \delta_{(T^{(u)}_i, u \Z_{ij})} \bigg|_{[0, 1] \times \EE_u}\,,
\]
in $[0, 1] \times \EE_u$, where
\begin{enumerate}
\item $\sum_i \delta_{T^{(u)}_i}$ is a homogeneous Poisson process on $[0, 1]$ with intensity $\theta u^{-\alpha}$;
\item $(\sum_j \delta_{\Z_{ij}})_i$ is an iid sequence of point processes in $\EE$, independent of $\sum_i \delta_{T^{(u)}_i}$, and with common distribution equal to the weak limit in \eqref{E:clusterprocess}.
\end{enumerate}
\end{thm}

In the setting of Theorem~\ref{T:pointprocess:complete}, the quantity $\theta$ in \eqref{E:candidatetheta} is indeed the extremal index of the sequence $(\norm{\X_t})$, that is, for all $u \in (0, \infty)$ and as $n \to \infty$,
\[
    \Pr (M_n \leq a_n u)
    = \{ \Pr (\norm{\X_1} \leq a_n u) \}^{n \theta} + o(1)
    \to e^{-\theta u^{-\alpha}}.
\]

%%%%%%%%%%%
\subsection{Moving averages with random coefficients}
\label{SS:main:MMA}

Consider the process $(\X_t)$ in \eqref{E:MMA}, a multivariate moving average process of finite order and with random coefficient matrices. Fix two arbitrary norms on $\RR^q$ and $\RR^d$, and on $\RR^{d \times q}$ consider the corresponding operator norm. Without any danger of confusion, all these norms are denoted by $\norm{\,\cdot\,}$.

Joint regular variation of the process $(\X_t)$ will be established under the following conditions:
\begin{itemize}
\item[(M1)] The law of $\xib_0$ is multivariate regularly varying of index $\alpha \in (0, \infty)$ and with spectral measure $\law(\Thetab)$ on $\Sb^{q-1}$.
\item[(M2)] There exists $\beta > \alpha$ such that $\E \norm{ \C_i(0)}^\beta < \infty$ for $i = 0, \ldots, m$.
\item[(M3)] For $\Thetab$ as in (M1) and independent of $\{\C_i(t)\}$, there exists $i = 0, \ldots, m$ such that $\Pr\{\norm{\C_i(0) \Thetab} > 0\} > 0$.
\end{itemize}

\begin{thm}
\label{T:MMA}
Let $(\X_t)$ be as in \eqref{E:MMA}. If (M1)--(M3) hold,
then $(\X_t)$ is jointly regularly varying of index $\alpha$. As $x \to \infty$,
\begin{equation}
\label{E:MMA:tailequiv}
    \frac{\Pr(\norm{\X_0} > x)}{\Pr(\norm{\xib_0} > x)}
    \to \sum_{i=0}^m \E \norm{\C_i(0) \Thetab}^\alpha =: c > 0,
\end{equation}
and for $s, t \in \ZZ$ with $s \leq t$ and bounded and continuous $f : (\RR^d)^{t-s+1} \to \RR$,
\begin{eqnarray}
\label{E:MMA:spectral}
    \lefteqn{
    \E \biggl[ f\biggl( \frac{\X_s}{\norm{\X_0}}, \ldots, \frac{\X_t}{\norm{\X_0}} \biggr)
    \, \bigg| \, \norm{\X_0} > x \biggr]
    } \\
    &\to& c^{-1}
    \sum_{i=0}^m \E \biggl[ f\biggl(
    \frac{\C_{i+s}(s) \Thetab}{\norm{\C_i(0) \Thetab}}, \ldots,
    \frac{\C_{i+t}(t) \Thetab}{\norm{\C_i(0) \Thetab}}
    \biggr) \norm{\C_i(0) \Thetab}^\alpha \biggr], \nonumber
\end{eqnarray}
where $\C_j(\,\cdot\,) = \0$ if $j < 0$ or $j > m$.
\end{thm}

In view of Corollary~\ref{cor:spectralprocess}, equation~\eqref{E:MMA:spectral} determines the law of the spectral process of $(\X_t)$ and thus also of the tail process itself. The candidate extremal index $\theta$ is computed in Remark~\ref{rem:MMA:theta}.

\section{Joint regular variation and the tail process}
\label{S:tailprocess}

The principal aim of this section is to give the proofs of the results of subsection~\ref{SS:main:tailprocess} on the connection between joint regular variation, the tail process and the spectral process. The section is closed by a number of additional properties of the spectral process in Theorems~\ref{T:tailprocess:bis} and  \ref{T:tailprocess:ter}.

For the reader's convenience, we first recall the concept of vague convergence of measures; see e.g.~\citet[chapter~15]{Kallenberg83} or \citet[section~3.4]{Resnick87}. Note that a subset $K$ of $\EE$ is compact if and only if it is closed as a subset of $[-\infty,\infty]^d$ and does not contain the origin; a Radon measure $\mu$ on $\EE$ is therefore a Borel measure such that $\mu(\EE_u) < \infty$ for every $u \in (0, \infty)$. A sequence of Radon measures $(\mu_n)$ on $\EE$ then converges vaguely to a Radon measure $\mu$ if $\int f \rmd \mu_n \to \int f \rmd \mu$ as $n \to \infty$ for every $f \in C_K^+(\EE)$. In the sequel, integrals will be denoted often in operator notation $\mu(f) = \int f \rmd \mu$.

In the course of the proof of Theorem~\ref{T:JRV}, the following equivalent characterization of multivariate regular variation as defined in \eqref{E:RV} will be needed. Recall that a measurable function $V : (0,\infty) \to (0,\infty)$ is regularly varying of index $\rho \in \RR$ if $V(xy) / V(x) \to y^\rho$ as $x \to \infty$ for all $y \in (0, \infty)$. A $d$-dimensional random vector $\X$ is then regularly varying of index $\alpha \in (0, \infty)$ if and only if there exists a regularly varying function $V$ of index $-\alpha$ and a nonzero Radon measure $\mu$ on $\EE$ such that, as $x \to \infty$,
\begin{equation}
\label{E:RV:mu}
   \frac{1}{V(x)} \Pr(x^{-1} \X \in \cdot\,) \vto \mu(\,\cdot\,)
\end{equation}
\citep[p.~69]{Resnick86}. The measure $\mu$ is homogeneous of order $-\alpha$; as a consequence, it does not put any mass on hyperplanes through infinity. A possible choice for the function $V$ in \eqref{E:RV:mu} is $V(x) = \Pr(\norm{\X} > x)$, in which case for all $u \in (0, \infty)$ and with $\Thetab$ as in \eqref{E:RV},
\begin{equation}
\label{E:homogeneous}
    \mu(\{ \x \mid \norm{\x} > u, \x / \norm{\x} \in \cdot \,\})
    = u^{-\alpha} \Pr(\Thetab \in \cdot \,).
\end{equation}

Part of the proof of Theorem~\ref{T:JRV} rests on the property that a certain class of functions $\mathcal{F} \subset C_K(\EE)$ is measure-determining, that is, two Radon measures $\mu$ and $\nu$ on $\EE$ coincide if and only if $\mu(f) = \nu(f)$ for every $f \in \mathcal{F}$. For the following lemma, fix $k, l \in \NN$ and identify $\EE^{k+l} = [-\infty, \infty]^{k+l} \setminus \{ \0 \}$ with $([-\infty, \infty]^k \times [-\infty, \infty]^l) \setminus \{ (\0, \0) \}$. Further, fix two arbitrary norms on $\RR^k$ and $\RR^l$, both of which are conveniently denoted by $\norm{\,\cdot\,}$.

\begin{lem}
\label{L:vague}
Every Radon measure $\mu$ on $\EE^{k+l}$ is uniquely determined by $\mu(f)$ with $f$ ranging over $\mathcal{F}_1 \cup \mathcal{F}_2$ where
\begin{eqnarray*}
    \mathcal{F}_1 &=& \{ f \in C_K(\EE^{k+l}) \mid \exists u \in (0, \infty) : \norm{\y_1} \leq u \Rightarrow f(\y_1, \y_2) = 0 \}, \\
    \mathcal{F}_2 &=& \{ f \in C_K(\EE^{k+l}) \mid f(\y_1, \y_2) = f(\0, \y_2) \}.
\end{eqnarray*}
\end{lem}

\begin{proof}
For $n \in \NN$, define $\phi_n(\y_1, \y_2) = \min \{\max (n \norm{\y_1} - 1, 0), 1\}$ where $(\y_1, \y_2) \in \EE^{k+l}$. Clearly $\phi_n \in \mathcal{F}_1$. Moreover, as $n \to \infty$, the sequence $\phi_n$ increases pointwise to the indicator function of $\mathbb{F}_{k,l} = \{ (\y_1, \y_2) \mid \y_1 \neq 0 \}$.

Let $\mu$ be a Radon measure on $\EE^{k+l}$ and let $g \in C_K(\EE^{k+l})$. We have to show that $\mu(g)$ is uniquely determined by the values of $\mu(f_i)$ with $f_i$ ranging over $\mathcal{F}_i$ and $i \in \{1, 2\}$. Define the function $f_2$ by $f_2(\y_1, \y_2) = g(\0, \y_2)$; clearly $f_2 \in \mathcal{F}_2$. The function $g$ can be decomposed as
\[
    g
    = (g - f_2) + f_2
    = (g - f_2) \1_{\mathbb{F}_{k,l}} + f_2
    = \lim_{n \to \infty} (g - f_2) \phi_n + f_2.
\]
By the dominated convergence theorem, $\mu\{(g-f_2) \phi_n\} + \mu(f_2) \to \mu(g)$ as $n \to \infty$. Since $(g - f_2) \phi_n \in \mathcal{F}_1$ for every $n$, the lemma follows.
\end{proof}

\begin{proof}[Proof of Theorem~\ref{T:JRV}]
\textit{(i) implies (iii).}
Without loss of generality, assume $s \leq 0 \leq t$. By assumption, the law of $(\X_s, \ldots, \X_t)$ is regularly varying at infinity of index $\alpha$. A possible choice for $V$ in \eqref{E:RV:mu} is $V(x) = \Pr (\max_{i = s, \ldots, t} \norm{\X_i} > x)$. The limit of $\Pr(\norm{\X_0} > x) / V(x)$ as $x \to \infty$ must exists and since $\Pr(\norm{\X_0} > x) \leq V(x) \leq (t-s+1) \Pr(\norm{\X_0} > x)$, this limit must be a finite, positive constant. Therefore, an alternative choice for $V$ in \eqref{E:RV:mu} is $V(x) = \Pr(\norm{\X_0} > x)$: there exists a non-trivial Radon measure $\mu_{s,t}$ on $\EE^{(t-s+1)d}$ such that as $x \to \infty$,
\begin{equation}
\label{E:must}
    \frac{1}{\Pr(\norm{\X_0} > x)}
    \Pr \{ (x^{-1} \X_s, \ldots, x^{-1} \X_t) \in \, \cdot \, \}
    \vto \mu_{s,t}(\,\cdot\,).
\end{equation}
By construction, the restriction of $\mu_{s,t}$ to the set $\{(\y_s, \ldots, \y_t) \mid \norm{\y_0} > 1\}$ is a probability measure, say $\nu_{s, t}$. Let $(\Y_s, \ldots, \Y_t)$ be a random vector with law $\nu_{s, t}$. For $f : (\RR^d)^{t-s+1} \to \RR$ bounded and continuous, as $x \to \infty$,
\begin{eqnarray*}
    \lefteqn{
    \E [ f(x^{-1} \X_s, \ldots, x^{-1} \X_t) \mid \norm{\X_0} > x ]
    } \\
    &\to& \int f(\y_s, \ldots, \y_t) \1(\norm{\y_0} > 1) \mu_{s,t}(\d \y)
    = \E [ f(\Y_s, \ldots, \Y_t) ].
\end{eqnarray*}
Here it was used that $\mu_{s,t}$ is homogeneous and therefore puts no mass on the set of $(\y_s, \ldots, \y_t)$ for which $\norm{\y_0} = 1$ as well as on the set of vectors with at least one infinite coordinate. The above display establishes the convergence in distribution stated in (iii). By Kolmogorov's extension theorem, there exists a single random process $(\Y_t)_{t \in \ZZ}$ such that for all $s, t \in \ZZ$, the distribution of $(\Y_s, \ldots, \Y_t)$ is $\nu_{s,t}$. The law of $\norm{\Y_0}$ follows from the fact that the function $x \mapsto \Pr(\norm{\X_0} > x)$ is regularly varying of index $-\alpha$.

\textit{(iii) implies (ii).}
Trivial.

\textit{(ii) implies (i).} For every $t \in \NN$, the vector $(\X_0, \ldots, \X_t)$ will be shown to be regularly varying of index $\alpha$ in the sense of \eqref{E:RV:mu} with $V(x) = \Pr(\norm{\X_0} > x)$. The proof is by induction on $t$. The case $t = 0$ is trivial. So let $t \geq 1$. For $x \in (0, \infty)$, define the finite measure $\mu_x$ on $\EE^{(t+1)d}$ by
\[
    \mu_x (\,\cdot\,) = \frac{1}{\Pr( \norm{\X_0} > x )} \Pr\{ (x^{-1} \X_0, \ldots, x^{-1} \X_t) \in \cdot \, \}.
\]
It has to be shown that $\mu_x \vto \mu$ as $x \to \infty$ for some Radon measure $\mu$ on $\EE^{(t+1)d}$. Since $\mu_x(\{(\y_0, \ldots, \y_t) \mid \norm{\y_0} > 1\}) = 1$, the vague limit, $\mu$, provided it exists, is certainly nonzero.

Vague convergence of $\mu_x$ will follow from the following two
statements:
\begin{itemize}
\item[(a)] the family $(\mu_x)_x$ is relatively compact in the vague topology;
\item[(b)] there can be at most one limit of $\mu_x$ as $x \to \infty$.
\end{itemize}

First, by \citet[Proposition~3.16]{Resnick87}, a necessary and sufficient condition for (a) is that $\sup_x \mu_x(B) < \infty$ for every Borel set $B$ with compact closure. For such $B$, there exists $u > 0$ such that $(\y_0, \ldots, \y_t) \in B$ implies $\norm{\y_i} > u$ for some $i \in \{0, \ldots, t\}$. But then, by stationarity,
\[
    \mu_x(B) \leq (t+1) \frac{\Pr(\norm{\X_0} > u x)}{\Pr(\norm{\X_0} > x)}.
\]
Since $\Pr(\norm{\X_0} > \cdot)$ is regularly varying of index $-\alpha$, (a) follows.

Second, to prove (b), it is sufficient to show that $\lim_{x \to
\infty} \mu_x(f)$ exists for every $f$ in a measure-determining subset $\mathcal{F}$
of $C_K(\EE^{(t+1)d})$. According to Lemma~\ref{L:vague} with $k = d$ and $l = td$, $\mathcal{F} = \mathcal{F}_1 \cup \mathcal{F}_2$ is such a set, where
\begin{eqnarray*}
    \mathcal{F}_1 &=& \{ f \mid \exists u \in (0,\infty) : f(\y_0, \ldots, \y_t) = 0 \mbox{ whenever } \norm{\y_0} \leq u \}, \\
    \mathcal{F}_2 &=& \{ f \mid f(\y_0, \y_1, \ldots, \y_t) = f(\0, \y_1, \ldots, \y_t) \}.
\end{eqnarray*}
On the one hand, if $f \in \mathcal{F}_1$ with $u$ as above, then
by (ii),
\begin{eqnarray*}
    \mu_x(f)
    &=& \frac{1}{\Pr(\norm{\X_0} > x)}
    \E [ f(x^{-1} \X_0, \ldots, x^{-1} \X_t) \1 ( \norm{x^{-1} \X_0} > u ) ] \\
    &=& \frac{\Pr(\norm{\X_0} > ux)}{\Pr(\norm{\X_0} > x)}
    \E [ f(x^{-1} \X_0, \ldots, x^{-1} \X_t) \mid \norm{\X_0} > ux ] \\
    &\to& u^{-\alpha} \E [ f(u\Y_0, \ldots, u\Y_t) ], \qquad x \to \infty.
\end{eqnarray*}
On the other hand, if $f \in \mathcal{F}_2$, then by stationarity,
\begin{eqnarray*}
    \mu_x(f)
    &=& \frac{1}{\Pr(\norm{\X_0} > x)}
    \E [ f(0, x^{-1} \X_1, \ldots, x^{-1} \X_t) ] \\
    &=& \frac{1}{\Pr(\norm{\X_0} > x)}
    \E [ f(0, x^{-1} \X_0, \ldots, x^{-1} \X_{t-1}) ],
\end{eqnarray*}
the limit of which as $x \to \infty$ exists by the induction
hypothesis.
\end{proof}

\begin{proof}[Proof of Theorem~\ref{T:tailprocess}]
\textit{(i)}
Let $s,t \in \ZZ$ with $s \leq 0 \leq t$. Consider the following subsets of $\EE^{(t-s+1)d}$:
\begin{eqnarray*}
    \EE_{s,t} &=& \{ (\y_s, \ldots, \y_t) \mid 0 < \norm{\y_0} < \infty \}, \\
    \Sb_{s,t} &=& \{ (\y_s, \ldots, \y_t) \mid \norm{\y_0} = 1 \}.
\end{eqnarray*}
Further, define a bijection $T : \EE_{s,t} \to (0, \infty) \times \Sb_{s, t}$ by
\[
    T(\y_s, \ldots, \y_t)
    = \biggl( \norm{\y_0},
    \biggl( \frac{\y_s}{\norm{\y_0}}, \ldots, \frac{\y_t}{\norm{\y_0}}
    \biggr) \biggr).
\]
Let $\mu_{s,t}$ be as in \eqref{E:must} and define the measure $\Phi_{s,t}$ on $\Sb_{s,t}$ by
\[
    \Phi_{s,t}(B) = \mu_{s,t} \bigl( T^{-1} ((1, \infty) \times B) \bigr)
\]
for Borel-measurable $B \subset \Sb_{s,t}$. Since the law of $(\Y_{s}, \ldots, \Y_t)$ is equal to the restriction of $\mu_{s,t}$ to $T^{-1}((1, \infty) \times \Sb_{s,t})$, the measure $\Phi_{s,t}$ is in fact equal to the law of $(\Y_s / \norm{\Y_0}, \ldots, \Y_t / \norm{\Y_0}) = (\Thetab_{s}, \ldots, \Thetab_t)$. Moreover, as $\mu_{s,t}$ is homogeneous of order $-\alpha$, for $u \in (0,\infty)$ and Borel sets $B \subset \Sb_{s, t}$,
\begin{equation}
\label{E:factorization}
    \mu_{s,t} \bigl( T^{-1} ((u, \infty) \times B) \bigr)
    = u^{-\alpha} \Phi_{s,t}(B).
\end{equation}
For $u \geq 1$, the left-hand side is equal to $\Pr\{ \norm{\Y_0} > u, (\Thetab_s, \ldots, \Thetab_t) \in B\}$, while the right-hand side is equal to $\Pr( \norm{\Y_0} > u ) \Pr\{ (\Thetab_s, \ldots, \Thetab_t) \in B \}$. As a consequence, $\norm{\Y_0}$ and $(\Thetab_{-s}, \ldots, \Thetab_t)$ are independent. Since $s$ and $t$ were arbitrary, (i) follows.

\textit{(ii)}
Let again $s,t \in \ZZ$ with $s \leq 0 \leq t$, and let $g : (\RR^d)^{t-s+1} \to \RR$ be bounded and continuous and with the property that $g(\y_s, \ldots, \y_t) = 0$ if $\y_0 = \0$. By stationarity and \eqref{E:must},
\begin{eqnarray*}
    \lefteqn{ \E [ g(\Y_{s-i}, \ldots, \Y_{t-i}) ] } \\
    &=& \lim_{x \to \infty} \frac{1}{\Pr(\norm{\X_0} > x)}
    \E [ g(x^{-1} \X_{s-i}, \ldots, x^{-1} \X_{t-i}) \1(\norm{\X_0} > x) ] \\
    &=& \lim_{x \to \infty} \frac{1}{\Pr(\norm{\X_0} > x)}
    \E [ g(x^{-1} \X_s, \ldots, x^{-1} \X_t) \1(\norm{\X_i} > x) ] \\
    &=& \int g(\y_s, \ldots, \y_t) \1(\norm{\y_i} > 1)
        \mu_{s \wedge i, t \vee i}(\d\y).
\end{eqnarray*}
By the assumed property of $g$, the region of integration can be restricted to $\{ \y \mid \y_0 \neq \0 \}$. By \eqref{E:factorization} applied to $\mu_{s \wedge i, t \vee i}$,
\begin{eqnarray}
\label{E:shift:10}
    \lefteqn{
    \E [ g(\Y_{s-i}, \ldots, \Y_{t-i}) ]
    } \\
    &=& \int_0^\infty
    \E [ g(r\Thetab_s, \ldots, r\Thetab_t) \1(r\norm{\Thetab_i} > 1) ]
    \d(-r^{-\alpha}). \nonumber
\end{eqnarray}
Let $f$ be as in (ii) of the theorem and define
\begin{equation}
\label{E:gandf}
    g(\y_{s \wedge i}, \ldots, \y_{t \vee i})
    = f \biggl( \frac{\y_s}{\norm{\y_i}}, \ldots,
    \frac{\y_t}{\norm{\y_i}} \biggr) (\norm{\y_i} \wedge 1).
\end{equation}
Since $\Thetab_j = \Y_j / \norm{\Y_0}$ and $\norm{\Y_0} > 1$,
\begin{eqnarray*}
    \E [ f(\Thetab_{s-i}, \ldots, \Thetab_{t-i}) ]
    &=& \E  \biggl[
            f   \biggl(
                    \frac{\Y_{s-i}}{\norm{\Y_{i-i}}}, \ldots,
                    \frac{\Y_{t-i}}{\norm{\Y_{i-i}}}
                \biggr)
            (\norm{\Y_{i-i}} \wedge 1)
            \biggr] \\
    &=& \E [ g(\Y_{(s \wedge i) - i}, \ldots, \Y_{(t \vee i) - i}) ].
\end{eqnarray*}
In combination with \eqref{E:shift:10} applied to this particular function $g$, it follows that $\E [ f(\Thetab_{s-i}, \ldots, \Thetab_{t-i}) ]$ is equal to
\[
    \int_0^\infty
    \E [ g(r\Thetab_{s \wedge i}, \ldots, r\Thetab_{t \vee i})
        \1(r\norm{\Thetab_i} > 1) ]
    \d(-r^{-\alpha}).
\]
By definition of $g$ in \eqref{E:gandf}, the above expression can be rewritten as
\[
    \int_0^\infty
    \E  \biggl[
        f   \biggl(
                \frac{\Thetab_s}{\norm{\Thetab_i}}, \ldots,
                \frac{\Thetab_t}{\norm{\Thetab_i}}
            \biggr)
        (\norm{r\Thetab_i} \wedge 1)
        \1(r\norm{\Thetab_i} > 1)
        \biggr]
    \d(-r^{-\alpha}).
\]
Apply Fubini's theorem and use the formula $\int_0^\infty \1(r\norm{\Thetab_i} > 1) \d(-r^{-\alpha}) = \norm{\Thetab_i}^\alpha$ to identify the above expression with the right-hand side of \eqref{E:shift}.
\end{proof}

\begin{proof}[Proof of Corollary~\ref{cor:spectralprocess}]
\textit{(i) implies (iii).} The stated convergence in distribution follows from Theorem~\ref{T:JRV}(iii) and the continuous mapping theorem. The connection between the tail process and the spectral process was already established in Theorem~\ref{T:tailprocess}(i).

\textit{(iii) implies (ii).} Trivial.

\textit{(ii) implies (i).} Let $Y$ be a random variable independent of $(\Thetab_t)_{t \in \NN}$ and for which $\Pr(Y > y) = y^{-\alpha}$ for $y \in [1, \infty)$. If we can show that as $x \to \infty$,
\begin{equation}
\label{E:spectralprocess}
    \law \biggl( \frac{\norm{\X_0}}{x}, \frac{\X_0}{\norm{\X_0}}, \ldots, \frac{\X_t}{\norm{\X_0}}
        \, \bigg| \, \norm{\X_0} > x \biggr)
    \dto \law( Y, \Thetab_0, \ldots, \Thetab_t )
\end{equation}
then joint regular variation of $(\X_t)$ will follow from the continuous mapping theorem and Theorem~\ref{T:JRV}(ii) with $\Y_t = Y \Thetab_t$. So let $y \in [1, \infty)$, $t \in \NN$ and $f : C((\RR^d)^{t+1}) \to \RR$ be bounded and continuous. We have
\begin{eqnarray*}
	\lefteqn{
	\E \biggl[ \1\biggl( \frac{\norm{\X_0}}{x} > y \biggr) 
	f \biggl( \frac{\X_0}{\norm{\X_0}}, \ldots, \frac{\X_t}{\norm{\X_0}} \biggr)
    \, \bigg| \, \norm{\X_0} > x 
	\biggr]
	} \\
	&=& \frac{\Pr(\norm{\X_0} > xy)}{\Pr(\norm{\X_0} > x)}
	\E \biggl[ 
	f \biggl( \frac{\X_0}{\norm{\X_0}}, \ldots, \frac{\X_t}{\norm{\X_0}} \biggr)
    \, \bigg| \, \norm{\X_0} > x y
	\biggr].
\end{eqnarray*}
By regular variation of $\Pr(\norm{\X_0} > \cdot\,)$ and by (ii), the right-hand side converges as $x \to \infty$ to $y^{-\alpha} \E [ f(\Thetab_0, \ldots, \Thetab_t) ] = \E [ \1(Y > y) f(\Thetab_0, \ldots, \Thetab_t) ]$. This finishes the proof of \eqref{E:spectralprocess} and therefore of the corollary.
\end{proof}

In the remainder of this section, let $(\X_t)_{t \in \ZZ}$ be a stationary process in $\RR^d$, jointly regularly varying of index $\alpha \in (0, \infty)$ and with tail process $(\Y_t)$ and spectral process $(\Thetab_t)$ with respect to a given norm $\norm{\,\cdot\,}$.

\begin{thm}
\label{T:tailprocess:bis}
For $t \in \ZZ$,
\[
    \lim_{\delta \downarrow 0} \lim_{x \to \infty}
    \Pr(\norm{\X_0} > \delta x \mid \norm{\X_t} > x)
    = \E \norm{\Thetab_t}^\alpha.
\]
\end{thm}

\begin{proof}
By Theorem~\ref{T:JRV} and stationarity, as $x \to \infty$,
\begin{eqnarray*}
    \lefteqn{
    \Pr(\norm{\X_0} > \delta x \mid \norm{\X_t} > x)
    } \\
    &=& \frac{\Pr(\norm{\X_0} > \delta x)}{\Pr(\norm{\X_t} > x)}
    \Pr(\norm{\X_t} > x \mid \norm{\X_0} > \delta x) 
    \to \delta^{-\alpha} \Pr( \delta \norm{\Y_t} > 1 ).
\end{eqnarray*}
From the spectral decomposition of $(\Y_t)$ and Fubini's theorem,
\[
%\begin{eqnarray*}
%    \lefteqn{
    \delta^{-\alpha} \Pr( \delta \norm{\Y_t} > 1 )
%    } \\
%    &=& \delta^{-\alpha} \int_1^\infty \Pr( \delta y \norm{\Thetab_t} > 1)
%    \d(-y^{-\alpha}) \\
    = \int_\delta^\infty \Pr( r \norm{\Thetab_t} > 1 ) \d(-r^{-\alpha})  
%    \\ 
%    &=& \E \biggl[
%	\int_\delta^\infty \1( r \norm{\Thetab_t} > 1) \d(-r^{-\alpha}) \biggr]
    = \E  \{ \min(\norm{\Thetab_t}, \delta^{-1}) \}^\alpha.
%\end{eqnarray*}
\]
Let $\delta \downarrow 0$ to conclude the proof.
\end{proof}

\begin{thm}
\label{T:tailprocess:ter}
Fix $t \in \ZZ$.

(i) For $h \in C(\Sb^{d-1})$,
\begin{equation}
\label{E:shift:0}
    \E [ h(\Thetab_t / \norm{\Thetab_t}) \norm{\Thetab_t}^\alpha ]
    = \E [ h(\Thetab_0) \1( \Thetab_{-t} \neq 0) ].
\end{equation}
and in particular $\E \norm{\Thetab_t}^\alpha = \Pr(\Thetab_{-t} \neq 0)$.

(ii) $\E \norm{\Thetab_t}^\alpha = 1$ if and only if for every $h \in C(\Sb^{d-1})$,
\[
    \E [ h(\Thetab_t / \norm{\Thetab_t}) \norm{\Thetab_t}^\alpha ]
    = \E [ h(\Thetab_0) ].
\]
\end{thm}

\begin{proof}
\textit{(i)}
Let $z \in (0, \infty)$ and define $f_z(\y) = h(\y/\norm{\y}) \min(\norm{\y}^\alpha, z)$ for $\y \in \RR^d \setminus \{ \0 \}$ while $f_z(\0) = 0$. The function $f_z$ is bounded and continuous. Hence, by Theorem~\ref{T:tailprocess}(ii) with $s = t = 0$, for $i \in \ZZ$,
\begin{eqnarray}
\label{E:shift:0:z}
    \E [ h(\Thetab_{-i}/\norm{\Thetab_{-i}}) \min(\norm{\Thetab_{-i}}^\alpha, z) ]
    &=& \E [ f_z(\Thetab_{-i}) ] \\
    &=& \E [ f_z(\Thetab_0 / \norm{\Thetab_i}) \norm{\Thetab_i}^\alpha ] \nonumber \\
    &=& \E [ h(\Thetab_0) \min( 1, z \norm{\Thetab_i}^\alpha) ]. \nonumber
\end{eqnarray}
The case $h \equiv 1$ yields $\E \min(\norm{\Thetab_{-i}}^\alpha, z) = \E \min( 1, z \norm{\Thetab_i}^\alpha)$, whence, by monotone convergence, $\E \norm{\Thetab_{-i}}^\alpha = \Pr(\Thetab_i \neq 0)$. Let $z \to \infty$ in \eqref{E:shift:0:z} to arrive at \eqref{E:shift:0} with $t$ replaced by $-i$.

\textit{(ii)}
The `if' part follows from $h \equiv 1$. The `only if' part follows from (i).
\end{proof}

% ===============================================================================
\section{Point processes}
\label{S:pointprocess}

Convergence of the time-space point process $N_n$ in \eqref{E:Nn} is often referred to in the literature as complete convergence. Complete convergence was claimed to hold in Theorem~\ref{T:pointprocess:complete}. The atoms of the limit point process of $N_n$ can be partitioned into independent and identically distributed clusters, the distribution of which can be written in terms of the tail process via Theorem~\ref{T:pointprocess:cluster}. Moreover, the Laplace functional of the cluster point process admits a succinct representation in terms of the forward tail process $(\Y_t)_{t \in \NN}$, see Theorem~\ref{T:Laplace}. All this is the content of subsection~\ref{SS:pointprocess:complete}.

Stripping the time or space coordinates from $N_n$ yields the derived point processes
\begin{eqnarray}
\label{E:Nnprime}
    N_n' &=& \sum_{i=1}^n \delta_{\X_i / a_n}, \\
\label{E:Nncirc}
	N_n^\circ &=& \sum_{i=1}^n \delta_{i/n} \1(\norm{\X _i} > a_n),
\end{eqnarray}
with state spaces $\EE$ and $[0, 1]$, respectively. Some known and new results on the weak limits of these processes are given in subsection~\ref{SS:pointprocess:derived}.

%%%%%%%%%%%
\subsection{Complete convergence and clusters of extremes}
\label{SS:pointprocess:complete}

For convenience, write $X_t = \norm{\X_t}$ and $Y_t = \norm{\Y_t}$ for $t \in \ZZ$ as well as $M_{i,j} = \max(X_i, \ldots, X_j)$ for $i, j \in \ZZ$ with $i \leq j$. Observe in particular that $M_r = M_{1,r}$ for $r \in \NN$.

\begin{proof}[Proof of Theorem~\ref{T:pointprocess:cluster}]
By condition~\ref{C:FiCl} and regular variation of the function $x \mapsto \Pr(X_0 > x)$, for all $u, v \in (0, \infty)$,
\be 
\label{AntiCl:v}
    \lim_{m\to\infty} \limsup_{n\to\infty}
    \Pr( M_{-r_n, -m} \vee M_{m, r_n} > a_n u v \mid X_0 > a_n u ) =
    0.
\ee
As a consequence, for every $v, \eps \in (0, \infty)$ there exists $m \in \NN$ such that for all $r \in \NN$ with $r \geq m$ we have $\Pr(\max_{m \leq |i| \leq r} Y_j > v) \leq \eps$. This proves that $\Pr( \lim_{|t| \to \infty} Y_t = 0 ) = 1$.

For $m, n \in \NN$, define
\begin{eqnarray}
\label{E:thetan}
    \theta_n &=& \frac{\Pr(M_{r_n} > a_n u)}{r_n \Pr(X_0 > a_n u)}, \\
\label{E:thetanm}
    \theta_{n,m} &=& \Pr(M_m \leq a_n u \mid X_0 > a_n u), \\
\nonumber
    \theta^Y_m &=& \Pr \biggl( \max_{i = 1, \ldots, m} Y_i \leq 1 \biggr),
\end{eqnarray}
and recall $\theta$ in \eqref{E:candidatetheta}. By \citet[section~2, p.~332]{Segers05}, $\liminf_{n \to \infty} \theta_n > 0$. Further, from the definition of the tail process, $\theta_{n,m} \to \theta^Y_m$ as $n \to \infty$, while by monotone convergence, $\theta^Y_m \to \theta$ as $m \to \infty$. Finally, by \citet[Theorem~3.1, eq.~(5)]{Segers05},
\[
    \lim_{m \to \infty} \limsup_{n \to \infty} |\theta_n - \theta_{n,m}| = 0.
\]
It follows that $\theta = \lim_{n \to \infty} \theta_n > 0$, as required. The proof of the identity $\theta = \Pr(E)$ is postponed until the end.

Consider now $f \in C_K^+(\EE)$. There exists $v \in (0, 1]$ such that $f(\x) = 0$ if $\norm{\x} \leq v$. For $s, t \in \ZZ$ such that $s \leq t$, write
\[
    c_n(s, t) = \exp \biggl\{ - \sum_{i=s}^t f(\X_i / (a_n u)) \biggr\}.
\]
We have
\[
    \E [ c_n(1, r_n) \mid M_{r_n} > a_n u ]
    = \frac{\E [ c_n(1, r_n) \1(M_{r_n} > a_n u) ]}{\Pr( M_{r_n} > a_n u )}.
\]
Split the event $\{M_{r_n} > a_n u\}$ according to the first time that the maximum is reached to get 
\begin{eqnarray}
\label{E:ClusterPeak:10}
	\lefteqn{
    \E [ c_n(1, r_n) \1(M_{r_n} > a_n u) ]
    } \\
    &=& \sum_{j=1}^{r_n} \E [ c_n(1, r_n) \1(a_n u \vee M_{j-1} \leq X_j = M_{r_n}) ]. \nonumber
\end{eqnarray} 
Fix $m \in \NN$ and let $n \in \NN$ be large enough so that $r_n \geq 2m+1$. For $j \in \NN$ such that $m+1 \leq j \leq r_n-m$, if $M_{j-m-1} \vee M_{j+m, r_n} \leq a_n u v$, then $c_n(1, r_n) = c_n(j-m, j+m-1)$ while $a_n u \vee M_{j-1} \leq X_j = M_{r_n}$ is equivalent to $a_n u \vee M_{j-m,j-1} \leq X_j = M_{j-m,j-m+1}$. Hence, for such $j$,
\begin{eqnarray}
\label{E:ClusterPeak:20}
    &&\bigl| \E [ c_n(1, r_n) \1(a_n u \vee M_{j-1} \leq X_j = M_{r_n}) ] \\
    && \quad \mbox{} - \E [ c_n(j-m, j+m-1) \1(a_n u \vee M_{j-m,j-1} \leq X_j = M_{j-m,j+m-1}) ] \bigr| \nonumber 
\end{eqnarray}
is bounded by
\begin{eqnarray*}
 \lefteqn{   \Pr( M_{j-m-1} \vee M_{j+m, r_n} > a_n u v, X_j > a_n u )}\\
    &\leq &\Pr(M_{-r_n,-m} \vee M_{m,r_n} > a_n u v, X_0 > a_n u),
\end{eqnarray*}
in view of stationarity, $0 < v \leq 1$, and $0 \leq c_n \leq 1$. By stationarity, the second expectation in \eqref{E:ClusterPeak:20} does not depend on $j$. Hence, in view of \eqref{E:ClusterPeak:10},
\begin{eqnarray*}
    &&\bigl| \E [ c_n(1, r_n) \1(M_{r_n} > a_n u) ] \\
    && \mbox{} \quad - r_n \E [ c_n(-m, m-1) \1(a_n u \vee M_{-m,-1} \leq X_0 = M_{-m,m-1}) ] \bigr| \\
    &\leq& 2 m \Pr(X_0 > a_n u) + r_n \Pr(M_{-r_n,-m} \vee M_{m,r_n} > a_n u v, X_0 > a_n u).
\end{eqnarray*}
Divide by $\Pr(M_{r_n} > a_n u)$ and recall $\theta_n$ in \eqref{E:thetan} to see that $\eps_{n,m}$ defined by
\begin{eqnarray*}
    \lefteqn{
    \bigl| \E [ c_n(1, r_n) \mid M_{r_n} > a_n u ]
    } \\ 
    && \mbox{} - \theta_n^{-1}
    \E [ c_n(-m, m-1)
    \1(a_n u \vee M_{-m,-1} \leq X_0 = M_{-m,m-1}) \mid X_0 > a_n u] \bigr|.
\end{eqnarray*}
is bounded by
\[
     \frac{r_n \Pr(X_j > a_n u)}{\Pr(M_{r_n} > a_n u)}
    \biggl( \frac{2 m}{r_n} 
    + \Pr(M_{-r_n,-m} \vee M_{m,r_n} > a_n u v \mid X_0 > a_n u) \biggr).
\]
From \eqref{AntiCl:v} and $\lim_n \theta_n = \theta > 0$, it follows that $\lim_m \limsup_n \eps_{n,m} = 0$. Therefore, by definition of the tail process, as $n \to \infty$,
\begin{eqnarray*}
    \lefteqn{
    \E [ c_n(1, r_n) \mid M_{r_n} > a_n u) ]
    } \\
    &\to& \theta^{-1} \E \biggl[ \exp \biggl( - \sum_{t \in \ZZ} f(\Y_i) \biggr)
    \1 \biggl( \sup_{t < 0} Y_t < Y_0 = \sup_{t \in \ZZ} Y_t \biggr) \biggr].
\end{eqnarray*}
The special case $f = 0$ yields $\theta = \Pr(E)$. This identity in combination with the previous display yields~\eqref{E:clusterprocess}.
\end{proof}

The description of the weak limit $\sum_j \delta_{\Z_j}$ in \eqref{E:clusterprocess} involves the distribution of the double sided tail process $(\Y_t)_{t\in\ZZ}$. In many cases the distribution of the forward tail process $(\Y_t)_{t \geq 0}$ is much more easily accessible than the one of the backward tail process, $(\Y_t)_{t \leq 0}$. This is the case for instance for Markov chains such as the random coefficient autoregressive process defined by $\X_t = \A_t \X_{t-1} + \B_t$ for $t \in \ZZ$, where $(\A_t, \B_t)$ are independent random elements in $\RR^{d \times d} \times \RR^d$ satisfying certain conditions \citep{Ke73}. Its forward tail process is a multiplicative random walk, $\Y_t = \A_t \cdots \A_1 \Y_0$ for positive integer $t$, but the backward tail process is to be constructed from the forward one via e.g.~Theorem~\ref{T:tailprocess}. It is interesting then that in general, the Laplace functional of $\sum_j \delta_{\Z_j}$ can also be described using the forward tail process only.

\begin{thm}
\label{T:Laplace}
The Laplace functional of $\sum_j \delta_{\Z_j}$ in \eqref{E:clusterprocess} is given by
\begin{eqnarray*}
    \E e^{-\sum_j f(\Z_j)}
    &=& \theta^{-1} \int_0^\infty 
    \E \biggl[ e^{-\sum_{i=0}^\infty f(y \Thetab_i)} \1\biggl(y \sup_{i \geq 0} \norm{\Thetab_i} > 1\biggr)\\
	&& \quad \qquad \qquad 
	- e^{-\sum_{i=1}^\infty f(y \Thetab_i)} 
	\1\biggl(y \sup_{i \geq 1} \norm{\Thetab_i} > 1\biggr) \biggr] 
	\rmd(-y^{-\alpha})
\end{eqnarray*}
for $f \in C_K^+(\EE)$. If additionally $f(\x) = 0$ whenever $\norm{\x} \leq 1$, then
\[
    \E e^{-\sum_j f(\Z_j)}
    = 1 - \theta^{-1} \E \Bigl[ e^{-\sum_{i=1}^\infty f(\Y_i)} - e^{-\sum_{i=0}^\infty
    f(\Y_i)} \Bigr]
\]
\end{thm}

\begin{proof}
Take $u \in (0, \infty)$ and $f\in C_K^+(\EE)$. There exists $v \in (0, 1]$ such that $f(\x) = 0$ whenever $\norm{\x} \leq v$. Hence
\begin{eqnarray*}
    \lefteqn{
    \E \Bigl[ e^{-\sum_{i=1}^{r_n} f(\X_i / (a_n u))} 
    			\,\Big|\, M_{r_n} > a_n u \Bigr]
    } \\
    &=& \frac{1}{\Pr( M_{r_n} > a_n u )}
      \E \Bigl[ e^{-\sum_{i=1}^{r_n} f(\X_i / (a_n u))} \1(M_{r_n} > a_n u) \Bigr] \\
    &=& \frac{\Pr( M_{r_n} > a_n u v )}{\Pr( M_{r_n} > a_n u )}
        \E \Bigl[ e^{-\sum_{i=1}^{r_n} f(\X_i / (a_n u))} \1(M_{r_n} > a_n u) 
        \,\Big|\, M_{r_n} > a_n u v \Bigr].
\end{eqnarray*}
By Theorem~\ref{T:pointprocess:cluster} and regular variation of $\Pr(\norm{\X_0} > \cdot\,)$, we have $\Pr( M_{r_n} > a_n u v ) / \Pr( M_{r_n} > a_n u ) \to v^{-\alpha}$ as $n \to \infty$. For $s, t \in \ZZ$ such that $s \leq t$, put
\[
    c_n(s,t) = e^{-\sum_{i=s}^{t} f(\X_i / (a_n u))} \1(M_{s,t} > a_n u).
\]
Note that $c_n(1,r_n) = c_n(\min I_n, \max I_n)$ with $I_n = \{i = 1, \ldots, r_n \mid \norm{\X_i} > a_n u v \}$. Hence, by an argument similar to the proof of Theorem~3.1 in \citet{Segers05} [see also the proof of Theorem~1 in \citet{Segers03}],
\[
    \lim_{m \to \infty} \limsup_{n \to \infty}
    | \E [ c_n(1, r_n) \mid M_{r_n} > a_n u v ] - \theta_{n,m}^{-1} A_{n,m} | = 0
\]
with $\theta_{n,m}$ as in \eqref{E:thetanm} and
\[
    A_{n,m} = \E [ c_n(0,m) - c_n(1,m) \1(M_m > a_n u v) 
    \mid \norm{\X_0} > a_n u v ].
\]
From the proof of Theorem~\ref{T:pointprocess:cluster}, $\lim_m \lim_n \theta_{n,m} = \theta$. Hence, by definition of the tail process
\begin{eqnarray*}
     \lim_{m \to \infty} \lim_{n \to \infty} A_{n,m}
    & = & \E \biggl[ e^{-\sum_{i=0}^\infty f(v\Y_i)} 
    \1\biggl(v \sup_{i \geq 0} \norm{\Y_i} > 1\biggr) \\
    & & \qquad - e^{-\sum_{i=1}^\infty f(v\Y_i)} 
    \1\biggl(v \sup_{i \geq 1} \norm{\Y_i} > 1\biggr) \biggr] = A.
\end{eqnarray*}
As a consequence, 
\[
    \E e^{-\sum_j f(\Z_j)}
    = \lim_{n \to \infty} \E \Bigl[ e^{-\sum_{i=1}^{r_n} f(\X_i / (a_n u))} 
    	\,\Big|\, M_{r_n} > a_n u \Bigr]
    = v^{-\alpha} \theta^{-1} A.
\]
From the spectral decomposition of the tail process $(\Y_t)$,
\begin{eqnarray*}
    v^{-\alpha} A
    &=& \int_v^\infty \E \biggl[ e^{-\sum_{i=0}^\infty f(y \Thetab_i)} 
    \1\biggl(y \sup_{i \geq 0} \norm{\Thetab_i} > 1\biggr) \\
    & & \qquad \qquad - e^{-\sum_{i=1}^\infty f(y \Thetab_i)} 
    \1\biggl(y \sup_{i \geq 1} \norm{\Thetab_i} > 1\biggr) \biggr] 
    \rmd(-y^{-\alpha}).
\end{eqnarray*}
As $\norm{\Thetab_0} = 1$ and $0 < v \leq 1$, the integrand is equal to zero for $0 < y \leq v$. Hence the domain of integration can be extended to $(0, \infty)$. The second formula follows from the first one upon noting that $\norm{\Thetab_0} = 1$ and $\theta = \Pr(\sup_{i \geq 1} \norm{\Y_i} \leq 1)$.
\end{proof}

\begin{proof}[Proof of Theorem~\ref{T:pointprocess:complete}]
Let $(\X_{k,j})_{j \in \NN}$, with $k \in \NN$, be independent copies of $(\X_j)_{j \in \NN}$, and define 
\[
	 \hat{N}_n 
	 = \sum_{k=1}^{k_n}  \hat{N}_{n,k}
	 \qquad \mbox{with} \qquad
	 \hat{N}_{n,k} = \sum_{j=1}^{r_n} 
	 \delta_{(k r_n/n, \X_{k,j}/a_n)}.
\]
By Condition~\ref{C:Aprimean}, the weak limits of $N_n$ and $\hat{N}_n$ must coincide. By \citet[Theorem 4.2]{Kallenberg83} it is enough to show that the Laplace functionals of $\hat{N}_n$ converge to those of $N^{(u)}$. Take $f \in C_K^+([0,1] \times \EE_u)$. It is convenient to adopt a convention that $f(t,\x)=0$ for all $(t,\x)\not\in \EE_u$. There exists $M \in (0, \infty)$ such that $0 \leq f(t, \x) \leq M \1(\norm{\x} > u)$. Hence as $n \to \infty$,
\begin{eqnarray*}
	1 
	\geq \E e^{-\hat{N}_{n,k} f} 
	&\geq& \E e^{-M\sum_{i=1}^{r_n} \1(\norm{\X_i} > a_nu)} \\
	&\geq& 1 - M r_n \Pr(\norm{\X_0} > a_n u) = 1 - o(1).
\end{eqnarray*}
In combination with the elementary bound $0 \leq - \log z - (1-z) \leq (1-z)^2/z$ for $z \in (0,1]$, it follows that as $n \to \infty$,
\[
	-\log \E e^{- \hat{N}_{n} f}
	= - \sum_{k=1}^{k_n} \log \E e^{- \hat{N}_{n,k} f}
	=  \sum_{k=1}^{k_n} (1 - \E e^{- \hat{N}_{n,k} f}) +o(1).
\]

By \eqref{E:runsblocks}, $k_n \Pr (M_{r_n} >a_n u ) \to \theta u^{-\alpha}$ for $u \in (0, \infty)$ and as $n \to \infty$. Hence
\begin{eqnarray*}
	\lefteqn{ 
	\sum_{k=1}^{k_n} (1 - \E  e^{- \hat{N}_{n,k} f}) 
	} \\
	&=&  k_n \Pr (M_{r_n}>a_n u ) 
	\sum_{k=1}^{k_n} \frac{1}{k_n}
	\E \biggl[ 1- e^{- \sum_{j=1}^{r_n} f(k r_n/n, \X_{j}/a_n) } \,\bigg|\,
	M_{r_n}>a_n u \biggr] \\
	&=&  \theta u^{-\alpha} \sum_{k=1}^{k_n} \frac{1}{k_n} 
	\E \biggl[ 1- e^{- \sum_{j=1}^{r_n} f(k r_n/n, \X_j/a_n)} 
	\,\bigg|\, M_{r_n}>a_n u \biggr] + o(1).
\end{eqnarray*}
Let $T_n$ be a random variable, uniformly distributed on $\{k r_n / n \mid k = 1, \ldots, k_n\}$ and independent of $(\X_j)$. By the previous display, as $n \to \infty$,
\[
	\sum_{k=1}^{k_n} (1 - \E  e^{- \hat{N}_{n,k} f}) 
	=  \theta u^{-\alpha}  \E \biggl[ 1- e^{- \sum_{j=1}^{r_n}
	f(T_n, u\X_j/(ua_n)) } \,\bigg|\, M_{r_n}>a_n u \biggr] + o(1).
\]
By \eqref{E:clusterprocess} and since $T_n$ converges in law to a random variable $T$ that is uniformly distributed on $(0,1)$, the expressions in the previous display converge as $n \to \infty$ to
\begin{equation}
\label{E:complete:limit}
	\theta u^{-\alpha}  
	\E \Bigl[ 1- e^{- \sum_j f(T, u \Z_j) } \Bigr]
	=  \int_0^1 \E \Bigl[ 1 - e^{-\sum_j f(t,u \Z_j) } \Bigr] 
	\theta u^{-\alpha} \rmd t.
\end{equation}
It remains to be shown that the right-hand side above equals $-\log \E e^{-N^{(u)}f}$  for $N^{(u)}$ as in the theorem.

Define $g(t) =\E \exp \{ - \sum_j f(t, u \Z_j) \}$ for $t \in [0, 1]$. Since $\sum_i \delta_{T_i^{(u)}}$ is independent of the iid sequence $(\sum_j \delta_{\Z_{ij}})_i$,
\[
	\E  e^{- N^{(u)}f} 
	= \E e^{- \sum_i \sum_j f(T^{(u)}_i, u\Z_{ij})}
	= \E e^{\sum_i \log g(T^{(u)}_i)}.
\]
The right-hand side is the Laplace functional of a homogeneous Poisson process on $[0,1]$ with intensity $\theta u^{-\alpha}$ evaluated in $- \log g$, which is equal to 
\[
	\exp \biggl( - \int_0^1 \{1 - g(t)\} \theta u^{-\alpha} \rmd t \biggr)
\]
see e.g. \citet[Lemma~5.1.12]{EKM}; note that $0 \leq g \leq 1$. By definition of $g$, the integral in the exponent is equal to the one in \eqref{E:complete:limit}.
\end{proof}

The fact that complete convergence in Theorem~\ref{T:pointprocess:complete} might hold was already mentioned without proof in \citet{DaHs95} with a reference to \citet{Mori}.

%%%%%%%%%%%
\subsection{Derived point processes}
\label{SS:pointprocess:derived}

Omitting the time component from the point processes $N_n$ in \eqref{E:Nn} yields the point processes $N_n'$ in \eqref{E:Nncirc}, living in the state space $\EE$. The limit behavior of $N_n'$ has been studied in \citet[section~2]{DaHs95} and \citet[section~2]{DaMi98} with the aim of determining asymptotics of sum-type functionals of $(\X_t)$ such as sample autocovariances and sample autocorrelations. Besides Condition~\ref{C:FiCl}, the following one is used in these papers.

\begin{cond}[$\mathcal{A}(a_n)$]
\label{C:Aan} There exists a positive integer sequence $(r_n)$
such that $r_n \to \infty$ and $r_n / n \to 0$ as $n \to \infty$
and such that for every $f \in C_K^+(\EE)$, denoting $k_n =
\lfloor n/r_n \rfloor$,
\[
    \E \exp \biggl( -\sum_{i=1}^nf(\X_i/a_n) \biggr) 
    - \biggl\{ \E \exp \biggl( -\sum_{i=1}^{r_n} f(\X_i/a_n) \biggr) \biggr\}^{k_n}
    \; \,\to\ \,  0\,.
\]
\end{cond}

Clearly $\mathcal{A}(a_n)$ is weaker than our condition $\mathcal{A'}(a_n)$, but both of them are satisfied for strongly mixing series. Note that in \citet{DaHs95} and \citet{DaMi98} step functions rather than continuous functions are used in the definition of $\mathcal{A}(a_n)$.

Using Theorem~\ref{T:pointprocess:cluster} the asymptotic behavior of $(N_n')$ can be described in somewhat more detail than what can be found in \citet{DaHs95} and \citet{DaMi98}. In particular, the limit of the point processes $(N_n')$ can be described via the tail process or spectral process of $(\X_t)$. The proof of the following theorem is similar to but simpler than the one of Theorem~\ref{T:pointprocess:complete} and is therefore omitted.

\begin{thm}
\label{T:pointprocess:space}
Under the assumptions of Theorem~\ref{T:pointprocess:cluster}, if also Condition~\ref{C:Aan} holds, then $N_n' \dto N'$ as $n \to \infty$ in $\EE$, where
\[
    N' = \sum_i \sum_j \delta_{P_i \Q_{ij}}
\]
consisting of the following ingredients:
\begin{enumerate}
\item a non-homogeneous Poisson process $\sum_i \delta_{P_i}$ on $(0, \infty)$ with intensity measure $\nu(\dy) = \theta \alpha y^{-\alpha-1} \dy$ for $y \in (0, \infty)$;
\item an iid sequence $(\sum_j \delta_{\Q_{ij}})_i$ of point processes in $\RR^d$, independent of $\sum_i \delta_{P_i}$, and with common law equal to the one of $\sum_j \delta_{\Z_j / M}$, where $M = \sup_j \norm{\Z_j}$.
\end{enumerate}
\end{thm}

\begin{rem}
\label{rem:clusters}
By the continuous mapping theorem, the common distribution of the point processes $\sum_j \delta_{\Q_{ij}}$ in the second item is equal to the weak limit as $n \to \infty$ in 
\[
    \law \Biggl( \sum_{i=1}^{r_n} \delta_{\X_i / M_{r_n}} 
    \,\Bigg|\, M_{r_n} > a_n u\Biggr)
    \dto \law \Biggl( \sum_{t \in \ZZ} \delta_{\Thetab_t} \,\Bigg|\, E \Biggr).
\]
Note that the event $E$ in \eqref{E:E} can be expressed in terms of the spectral process as well.
\end{rem}

Stripping the space component from the processes $N_n$ in \eqref{E:Nn} yields the processes $N_n^\circ$ in \eqref{E:Nncirc} with state space $[0,1]$. It is well-known that under appropriate mixing conditions, the processes $N_n^\circ$ converge weakly to a compound Poisson process \citep{LeRo88, HsHuLe89}. The distribution of the cluster sizes has been derived for several special Markovian models \citep{deHa89, KlPe}.

From Theorem~\ref{T:pointprocess:complete} with $u=1$ and \citet[Theorem~4.2]{Kallenberg83}, the limit behavior of $N_n^\circ$ follows at once. The distribution of the cluster sizes can be described in terms of the random variable
\begin{equation}
\label{E:nu}
	\nu = \sum_{i=1}^\infty \1(\norm{\Y_i} > 1).
\end{equation}
Note that $\theta = \Pr(\nu = 0)$. 
 
\begin{cor}
\label{cor:clustersizes} 
Under the assumptions of Theorem~\ref{T:pointprocess:complete}, $N_n^\circ \to N^\circ$ as $n \to \infty$ in $[0, 1]$, where $N^\circ = \sum_i \kappa_k \delta_{T_i}$ is a compound Poisson process consisting of: a homogenous Poisson process $\sum_i \delta_{T_i}$ on $[0,1]$ with intensity $\theta$; an iid sequence $(\kappa_i)$ of positive integer valued random variables, independent of $\sum_i \delta_{T_i}$, and with common law equal to the one of $\kappa = \sum_j \1(\norm{\Z_j} > 1)$. Moreover, for $s \in [0, \infty)$ and for integer $k \geq 1$,
\begin{eqnarray*}
	\E e^{-s \kappa} &=& 1 - (1 - e^{-s}) \theta^{-1} \E e^{-s \nu}, \\
	\Pr(\kappa = k) &=& \theta^{-1} \{ \Pr(\nu = k-1) - \Pr(\nu = k) \}.
\end{eqnarray*}
\end{cor}

\begin{proof}
Only the last statement requires some explanation: The first formula is a consequence of Theorem~4.1 with $f(\x) = s \1(\norm{\x} > 1)$; note that $\Pr(\norm{\Y_i} = 1) = 0$ for all $i \in \ZZ$ and that $\sum_{i=0}^\infty \1(\norm{\Y_i} > 1) = 1 + \nu$. The second formula follows from the first one by properties of probability generating functions together with Leibniz' product rule.
\end{proof}

\begin{rem}[univariate processes]
\label{rem:univariate}
Let $(X_t)_{t \in \ZZ}$ be a stationary univariate time series, jointly regularly varying of index $\alpha \in (0, \infty)$ and with tail process $(Y_t)$ and spectral process $(\Theta_t)$, the norm being of course the absolute value $|\,\cdot\,|$. By construction, the random variable $\Theta_0$ takes values on the zero-dimensional unit sphere $\Sb^0 = \{ - 1, 1 \}$, and as $x \to \infty$,
\[
	\frac{\Pr(X_0 > x)}{\Pr(|X_0| > x)} \to \Pr(Y_0 > 1) = \Pr(\Theta_0 = 1) =: p.
\] 
Denote $X_t^+ = \max(X_t, 0)$, $Y_t^+ = \max(Y_t, 0)$, and $\Theta_t^+ = \max(\Theta_t, 0)$. If $p > 0$, then the process $(X_t^+)$ is jointly regularly varying of index $\alpha \in (0, \infty)$ as well; the law of its tail process is equal to the conditional law of $(Y_t^+)$ given $Y_0 > 1$, while the law of its spectral process is equal to the conditional law of $(\Theta_t^+)$ given $\Theta_0 = 1$. By \eqref{E:theta:spectral}, the candidate extremal index of $(|X_t^+|) = (X_t^+)$ and therefore also of $(X_t)$ itself is given by
\begin{eqnarray*}
	\lefteqn{
	\lim_{t \to \infty} \lim_{x \to \infty} 
	\Pr \biggl( \max_{i = 1, \ldots, t} X_i \leq x
	\, \bigg| \, X_0 > x \biggr)
	} \\
	&=& \Pr \biggl( \sup_{i \geq 1} Y_i < 1 \, \bigg| \, Y_0 > 1 \biggr) 
	% \\ &=& 
	= \E \biggl[ \sup_{i \geq 0} (\Theta_i^+)^\alpha 
	- \sup_{i \geq 1} (\Theta_i^+)^\alpha 
	\, \bigg| \, \Theta_0 = 1 \biggr].
\end{eqnarray*}
Cluster size probabilities of $(X_t)$ are to be derived via Corollary~\ref{cor:clustersizes} from the law of $\sum_{i \geq 1} \1(Y_i > 1)$ conditionally on $Y_0 > 1$.
\end{rem}

% ===============================================================================
\section{Moving averages with random coefficients}
\label{S:mma}

In this section, Theorem~\ref{T:MMA} is proven by means of a version of Breiman's (1965) lemma adapted to regularly varying processes. 
%The name of this lemma is a reference to \citet{Breiman}, who shows that multiplying a random variable with a regularly varying law by an independent random variable with a lighter tail produces a random variable whose law is again regularly varying and of the same index as the first one. 
The following version for multivariate regular variation appears as Proposition~A.1 in \citet{BDM02b}.

\begin{lem}
\label{L:Breiman}
Let $\Z$ be a $k$-dimensional random column vector and let $\A$ be a $d \times k$ random matrix, independent of $\Z$. Assume that $\Z$ is multivariate regularly varying of index $\alpha \in (0, \infty)$, i.e.\ there exist a regularly varying function $V$ of index $-\alpha$ and a nonzero Radon measure $\mu$ on $\EE^k$ such that as $x \to \infty$,
\[
	\frac{1}{V(x)} \Pr(x^{-1} \Z \in \cdot \,) \vto \mu(\,\cdot\,).
\]
If $\E \norm{\A}^\beta < \infty$ for some $\beta > \alpha$, then in $\EE^d$, as $x \to \infty$,
\[
	\frac{1}{V(x)} \Pr(x^{-1} \A \Z \in \cdot \,) 
	\vto \E[\mu \circ \A^{-1}(\,\cdot\,)].
\]
\end{lem} 

%The limit measure is to be read as $\E[\mu \circ \A^{-1}(\,\cdot\,)] = \E[\mu(\{ \x \mid \A\x \in \,\cdot\,\})]$. 
If $V(x) = \Pr(\norm{\Z} > x)$, then there exists a random vector $\Thetab$ on $\Sb^{k-1}$ such that $\mu(f) = \E [ \int_0^\infty f(u\Thetab) \rmd(-u^{-\alpha})]$ for $\mu$-integrable functions $f$. Therefore, for independent copies of $\A$ and $\Thetab$ and as $x \to \infty$,
\begin{equation}
\label{E:tailequiv}
	\frac{\Pr(\norm{\A \Z} > x)}{\Pr(\norm{\Z} > x)}
	\to \E[ \mu( \{ \x \mid \norm{\A \x} > 1 \}) ] \\
	%&=& \E \biggl[ \int_0^\infty 
	%\1 (\norm{u \A \thetab} > 1) \rmd(-u^{-\alpha}) \biggr] 
	= \E \norm{\A \Thetab}^\alpha. %\nonumber
\end{equation}
Note that $\E \norm{\A \Thetab}^\alpha \leq \E \norm{\A}^\alpha < \infty$. If additionally $\Pr(\norm{\A\Thetab} > 0) > 0$, then also $\E \norm{\A \Thetab}^\alpha > 0$, so that $\norm{\A \Z}$ and $\norm{\Z}$ are tail equivalent. The following result provides a version of Breiman's lemma for processes.

\begin{lem}
\label{L:Breiman:proc}
Let $(\Z_t)_{t \in \ZZ}$ be a stationary sequence of random column vectors in $\RR^k$ and let $(\A_t)_{t \in \ZZ}$ be a stationary sequence of random $d \times k$ matrices, independent of $(\Z_t)$. Assume $(\Z_t)$ is regularly varying of index $\alpha \in (0, \infty)$ and spectral process $(\Thetab_t)$. If, for independent copies of $(\A_t)$ and $(\Thetab_t)$,
\begin{itemize}
\item[(a)] $\E \norm{\A_0}^\beta < \infty$ for some $\beta > \alpha$,
\item[(b)] $\Pr(\norm{\A_0 \Thetab_0} > 0) > 0$,
\end{itemize}
then $(\A_t \Z_t)$ is regularly varying of index $\alpha$ as well, and for $s, t \in \ZZ$ with $s \leq t$ and $f : (\RR^d)^{t-s+1} \to \RR$ bounded and continuous, as $x \to \infty$,
\begin{eqnarray}
\label{E:Breiman:spectral}
	\lefteqn{
	\E \biggl[f \biggl(\frac{\A_s \Z_s}{\norm{\A_0 \Z_0}}, \ldots, \frac{\A_t \Z_t}{\norm{\A_0 Z_0}} \biggr) \, \bigg| \, \norm{\A_0 \Z_0} > x \biggr]
	} \\
	&\to& \frac{1}{\E \norm{\A_0 \Thetab_0}^\alpha}
	\E \biggl[f \biggl(\frac{\A_s \Thetab_s}{\norm{\A_0 \Thetab_0}}, \ldots,
	\frac{\A_t \Thetab_t}{\norm{\A_0 \Thetab_0}}\biggr) 
	\norm{\A_0 \Thetab_0}^\alpha \biggr]. \nonumber
\end{eqnarray}
\end{lem}

\begin{proof}
Let $h : (\RR^d)^{t-s+1} \to \RR$ be bounded and continuous. In view of eqs.~\eqref{E:factorization} and \eqref{E:tailequiv} as well as Lemma~\ref{L:Breiman}, as $x \to \infty$,
\begin{eqnarray*}
	\lefteqn{
	\E[h(x^{-1} \A_s \Z_s, \ldots, x^{-1} \A_t \Z_t) \mid \norm{\A_0 \Z_0} > x]
	} \\
	&\to& \frac{1}{\E\norm{\A_0 \Thetab_0}^\alpha} 
	\int_0^\infty
	\E[h(u \A_s \Thetab_s, \ldots, u \A_t \Thetab_t) \1(u\norm{\A_0 \Thetab_0} > 1)] 
	\rmd (-u^{-\alpha}).
\end{eqnarray*}
Apply this relation to the function $h(\x_s, \ldots, \x_t) = f(\x_s / \norm{\x_0}, \ldots, \x_t / \norm{\x_0})$ to see that as $x \to \infty$, the left-hand side of \eqref{E:Breiman:spectral} converges to
\[
	\frac{1}{\E \norm{\A_0 \Thetab_0}^\alpha}
	\int_0^\infty \E \biggl[f \biggl(\frac{\A_s \Thetab_s}{\norm{\A_0 \Thetab_0}},
	\ldots, \frac{\A_t \Thetab_t}{\norm{\A_0 \Thetab_0}}\biggr) 
	\1(u\norm{\A_0 \Thetab_0} > 1) \biggr] \rmd(-u^{-\alpha}).
\]
By Fubini's theorem, this is equal to the right-hand side of \eqref{E:Breiman:spectral}.
\end{proof}

\begin{rem}[candidate extremal index]
In the setting of Lemma~\ref{L:Breiman:proc}, the candidate extremal index of $(\A_t \Z_t)$ in \eqref{E:candidatetheta}--\eqref{E:theta:spectral} is equal to
\begin{eqnarray}
\label{E:theta:Breiman}
	\lefteqn{
	\lim_{t \to \infty} \lim_{x \to \infty} 
	\Pr \biggl( \max_{i = 1, \ldots, t} \norm{\A_i \Z_i} \leq x
	\, \bigg| \, \norm{\A_0 \Z_0} > x \biggr)
	}  \nonumber \\
	&=& \frac{\E [ \sup_{i \geq 0} \norm{\A_i \Thetab_i}^\alpha
	- \sup_{i \geq 1} \norm{\A_i \Thetab_i}^\alpha ]}
	{\E \norm{\A_0 \Thetab_0}^\alpha}.
\end{eqnarray}
\end{rem}

\begin{rem}[linear combinations]
\label{rem:linear}
Let $(\X_t)$ be a stationary sequence of random $d$-dimensional column vectors, regularly varying of index $\alpha \in (0, \infty)$ and with spectral process $(\Thetab_t)$. Let $\a$ be a nonzero $d$-dimensional column vector. By \eqref{E:tailequiv}, as $x \to \infty$,
\[
	\frac{\Pr(|\a'\X_0| > x)}{\Pr(\norm{\X_0} > x)} \to \E |\a'\Thetab_0|^\alpha.
\]
If $\a'\Thetab_0$ is not degenerate at zero, then $\E |\a'\Thetab_0| > 0$, and by Lemma~\ref{L:Breiman:proc}, the univariate process $(\a'\X_t)$ is jointly regularly varying of index $\alpha$, the law of its spectral process $(\Theta_t^{\a})$ being given by
\[
	\E [ f (\Theta_s^{\a}, \ldots, \Theta_t^{\a}) ]
	= \frac{1}{\E|\a'\Thetab_0|^\alpha}
	\E \biggl[ f \biggl( \frac{\a'\Thetab_s}{|\a'\Thetab_0|}, \ldots, \frac{\a'\Thetab_t}{|\a'\Thetab_0|} \biggr) |\a'\Thetab_0|^\alpha \biggr]
\]
for integer $s, t \in \ZZ$ with $s \leq t$ and for integrable $f : \RR^{t-s+1} \to \RR$. By \eqref{E:theta:Breiman},
\begin{eqnarray}
\label{E:linear:theta:twosided}
	\lefteqn{
	\lim_{t \to \infty} \lim_{x \to \infty} 
	\Pr \biggl( \max_{i = 1, \ldots, t} |\a' \X_i| \leq x
	\, \bigg| \, |\a' \X_0| > x \biggr)
	} \nonumber \\
	&=& \frac{\E [ 	\sup_{i \geq 0} |\a'\Thetab_i|^\alpha 
					- \sup_{i \geq 1} |\a'\Thetab_i|^\alpha]}
	{\E |\a'\Thetab_0|^\alpha}. 
\end{eqnarray}
Similarly, by Remark~\ref{rem:univariate}, if $\Pr(\a'\Thetab_0 > 0) > 0$, writing $(z)_+^\alpha = \{\max(z,0)\}^\alpha$,
%the candidate extremal index of $(\a'\X_t)$ is given by
\begin{eqnarray}
\label{E:linear:theta:onesided}
	\lefteqn{
	\lim_{t \to \infty} \lim_{x \to \infty} 
	\Pr \biggl( \max_{i = 1, \ldots, t} \a' \X_i \leq x
	\, \bigg| \, \a' \X_0 > x \biggr)
	} \nonumber \\
	&=& \frac{\E [ \sup_{i \geq 0} (\a'\Thetab_i)_+^\alpha 
	- \sup_{i \geq 1} (\a'\Thetab_i)_+^\alpha ]}
	{\E (\a'\Thetab_0)^\alpha_+}.	
\end{eqnarray}
If Conditions~\ref{C:FiCl} and \ref{C:Aprimean} (or weaker versions tailored to $\a$) hold, then \eqref{E:linear:theta:twosided} and \eqref{E:linear:theta:onesided} are the extremal indices of $(|\a'\X_t|)$ and $(\a'\X_t)$, respectively.
\end{rem}

\begin{proof}[Proof of Theorem~\ref{T:MMA}]
We have $\X_t = \sum_{i=0}^m \C_i(t) \xib_{t-i} = \A_t \Z_t$ with $\A_t = (\C_0(t), \ldots, \C_m(t))$ a random matrix of dimension $d \times k$ where $k = (m+1)q$ and $\Z_t = (\xib_t', \ldots, \xib_{t-m}')'$ a random column vector of length $k$. The processes $(\A_t)$ and $(\Z_t)$ are stationary and independent of each other.

As the random vectors $\xib_t$ are mutually independent, it is straightforward to determine the tail process of $(\Z_t)$. First we specify the norms used in the sequel. On $\RR^k \cong (\RR^q)^{m+1}$, consider the norm $\norm{(\x_0, \ldots, \x_m)} = \max_{i=0,\ldots,m} \norm{\x_i}$ constructed from the chosen norm on $\RR^q$. The corresponding operator norm on $\RR^{d \times k} \cong (\RR^{d \times q})^{m+1}$ is given by $\norm{(\cc_0, \ldots, \cc_m)} = \max_{i = 0, \ldots, m} \norm{\cc_i}$, constructed in the same way from the operator norm on $\RR^{d \times q}$. Further, for $i = 0, \ldots, m$, let $\bfe_i$ be the $k \times q$ matrix
\[
    \bfe_i = (\0, \ldots, \0, I_q, \0, \ldots, \0)',
\]
where $\0$ and $I_q$ represent the $q \times q$ zero and identity matrices, respectively, $I_q$ appearing at position $i$. For $i \in \ZZ$ such that $i < 0$ or $i > m$, let $\bfe_i$ be the $k \times q$ zero matrix.

Let $Y$ be a random variable with survival function $\Pr(Y > y) = y^{-\alpha}$ for $y \in [1, \infty)$ and independent of $\Thetab$ in (M1). Put $\Y = Y \Thetab$. Assumption (M1) entails that $\law(x^{-1} \xib_0 \mid \norm{\xib_0} > x) \dto \law(\Y)$ as $x \to \infty$. Let $M$ be uniformly distributed on $\{0, 1, \ldots, m\}$ and independent of $Y$,  $\Thetab$, and $\{\C_i(t)\}$. Since the sequence $(\xib_t)$ is iid, $\law(x^{-1} \Z_0 \mid \norm{\Z_0} > x) \dto \law(\bfe_M \Y)$ as $x \to \infty$. Note that $\bfe_i \Y$ is a column vector of length $k = (m+1)q$ of which all entries are equal to zero except for those from position $iq+1$ to $iq+q$, which coincide with the entries of $\Y$. Put $\Y_t = \bfe_{M+t} \Y$ for $t \in \ZZ$. Then for $s, t \in \ZZ$ with $s \leq t$ and as $x \to \infty$,
\[
    \law(x^{-1} \Z_s, \ldots, x^{-1} \Z_t \mid \norm{\Z_0} > x)
    \dto \law(\Y_s, \ldots, \Y_t).
\]
Observe that $\Y_t = \0$ for $t \in \ZZ$ such that $|t| > m$, which is intuitively obvious from the construction of $(\Z_t)$. Since $\norm{\bfe_i \Thetab} = \norm{\Thetab} = 1$ for $i = 0, \ldots, m$, the spectral process of $(\Z_t)$ is simply $\Thetab_t = \bfe_{M+t} \Thetab$ for $t \in \ZZ$.

Having established joint regular variation of $(\Z_t)$, we only need to apply Lemmas~\ref{L:Breiman} and \ref{L:Breiman:proc}. Conditions~(a) and (b) of Lemma~\ref{L:Breiman:proc} follow from conditions (M2) and (M3) of the theorem, respectively. Note that for $t \in \ZZ$,
\begin{equation}
\label{E:ATheta}
	\A_t \Thetab_t 
	= (\C_0(t), \ldots, \C_m(t)) \bfe_{M+t} \Thetab 
	= \C_{M+t}(t) \Thetab,
\end{equation}
where $\C_i(t) = \0$ if $i < 0$ or $i > m$. By \eqref{E:tailequiv}, as $x \to \infty$,
\[
	\frac{\Pr(\norm{\X_0} > x)}{\Pr(\norm{\Z_0} > x)} 
	\to \E \norm{\A_0 \Thetab_0}^\alpha 
	= \frac{1}{m+1} \sum_{i=0}^m \E \norm{\C_i(0) \Thetab}^\alpha.
\]
Equation~\eqref{E:MMA:tailequiv} now follows from $\Pr(\norm{\Z_0} > x) \sim (m+1) \Pr(\norm{\xib_0} > x)$ as $x \to \infty$. Further, by \eqref{E:Breiman:spectral} and \eqref{E:ATheta}, the left-hand side of \eqref{E:MMA:spectral} converges to
\[	
	\frac{1}{\E \norm{\C_M(0) \Thetab}^\alpha}
	\E \biggl[ f\biggl( 
	\frac{\C_{M+s}(s) \Thetab}{\norm{\C_M(0) \Thetab}}, \ldots,
	\frac{\C_{M+t}(t) \Thetab}{\norm{\C_M(0) \Thetab}} 
	\biggr) \norm{\C_M(0) \Thetab}^\alpha \biggr].
\]
Condition on the value of $M$ to arrive at the right-hand side of \eqref{E:MMA:spectral}.
\end{proof}

\begin{rem}[candidate extremal index]
\label{rem:MMA:theta}
By \eqref{E:theta:Breiman} and the proof of Theorem~\ref{T:MMA}, for the moving average process $(\X_t)$ in \eqref{E:MMA}, the candidate extremal index in \eqref{E:candidatetheta} is equal to
\[
	\theta
    = \frac{\sum_{i=0}^m \E [ \sup_{t \geq 0} \norm{\C_{i+t}(t) \Thetab}^\alpha 
        - \sup_{t \geq 1} \norm{\C_{i+t}(t) \Thetab}^\alpha ]}
        {\sum_{i=0}^m \E \norm{\C_i(0) \Thetab}^\alpha}.
\]
Since $(\C_0(t), \ldots, \C_m(t))$ is stationary when indexed over $t \in \ZZ$, the numerator on the right-hand side is a telescoping sum, whence
\begin{equation}
\label{E:MMA:theta}
	\theta = \frac{\E [ \max_{i = 0, \ldots, m} \norm{\C_i(i) \Thetab}^\alpha ]}
        {\sum_{i=0}^m \E \norm{\C_i(0) \Thetab}^\alpha}.
\end{equation}
\end{rem}

\begin{rem}[finite-cluster condition]
\label{rem:MMA:FiCl}
The moving average $(\X_t)$ in \eqref{E:MMA} satisfies the finite-cluster condition~\ref{C:FiCl} under the following additional moment restriction on $\norm{\C_i(0)}$:
\begin{itemize}
\item[(M2')] For all $\gamma \in (0, 2\alpha)$ and all $i \in \{0, 1, \ldots, m\}$, we have $\E\norm{\C_i(0)}^\gamma < \infty$.
\end{itemize}
More precisely, under the assumptions of Theorem~\ref{T:MMA} with (M2) replaced by (M2'), Condition~\ref{C:FiCl} holds for every integer sequence $r_n \to \infty$ for which there exists $\eps \in (0, 1)$ such that $r_n = O(n^{1-\eps})$ as $n \to \infty$. The proof is straightforward and can be obtained from the authors. Of course, if the process $\{\C_i(t)\}$ is row-wise independent, then the moving average $(\X_t)$ is itself $(m+1)$-dependent, so that Conditions~\ref{C:FiCl} and \ref{C:Aprimean} both hold and $\theta$ in \eqref{E:MMA:theta} is the extremal index of $(\norm{\X_t})$.
\end{rem}

%%%%%%%%%%%

\end{document}